\documentclass[12pt,twoside]{article}
\usepackage[english]{babel}
\usepackage[latin1]{inputenc}
\usepackage{amsmath}
\usepackage{graphicx}
\usepackage{times,amssymb,amscd}

\numberwithin{equation}{section}

\newcommand{\bC}{\mathbf{C}}

\newcommand{\bG}{\mathbf{G}}
\newcommand{\bH}{\mathbf{H}}

\newcommand{\bL}{\mathbf{L}}

\newcommand{\bR}{\mathbf{R}}
\newcommand{\bS}{\mathbf{S}}

\newcommand{\br}{\mathbf{r}}

\newcommand{\bT}{\mathbf{T}}

\newcommand{\bt}{\mathbf{t}}

\newcommand{\cP}{\mathcal{P}}
\newcommand{\cS}{\mathcal{S}}

\newcommand{\cC}{\mathcal{C}}
\newcommand{\cH}{\mathcal{H}}

\newcommand{\EUC}{\mathbf E^3}
\newcommand{\SPH}{\bS^3}
\newcommand{\HYP}{\bH^3}
\newcommand{\SXR}{\bS^2\!\times\!\bR}
\newcommand{\HXR}{\bH^2\!\times\!\bR}
\newcommand{\SLR}{\widetilde{\bS\bL_2\bR}}
\newcommand{\NIL}{\mathbf{Nil}}
\newcommand{\SOL}{\mathbf{Sol}}

\begin{document}
\pagestyle{myheadings}
\markboth{\centerline{Jen\H o Szirmai}}
{Menelaus' and Ceva's theorems $\dots$}
\title
{Menelaus' and Ceva's theorems for translation triangles in Thurston geometries
}

\author{Jen\H o Szirmai \\
\normalsize Department of Algebra and Geometry, Institute of Mathematics,\\
\normalsize Budapest University of Technology and Economics, \\
\normalsize M\"uegyetem rkp. 3., H-1111 Budapest, Hungary \\
\normalsize szirmai@math.bme.hu
\date{\normalsize{\today}}}
\maketitle
\begin{abstract}
After having investigated and defined the ``surface of a translation-like triangle" in each
non-constant curvature Thurston geometry \cite{Cs-Sz25},
we generalize the famous Menelaus' and Ceva's theorems for translation triangles in the mentioned spaces.

The described method makes it possible to transfer further classical Euclidean theorems and notions to Thurston
geometries with non-constant curvature.
In our work we will use the projective models of Thurston geometries described
by E. Moln\'ar in \cite{M97}.

\end{abstract}
\newtheorem{Theorem}{Theorem}[section]
\newtheorem{corollary}[Theorem]{Corollary}
\newtheorem{lemma}[Theorem]{Lemma}
\newtheorem{exmple}[Theorem]{Example}
\newtheorem{definition}[Theorem]{Definition}
\newtheorem{rmrk}[Theorem]{Remark}
\newtheorem{proposition}[Theorem]{Proposition}
\newenvironment{remark}{\begin{rmrk}\normalfont}{\end{rmrk}}
\newenvironment{example}{\begin{exmple}\normalfont}{\end{exmple}}
\newenvironment{acknowledgement}{Acknowledgement}

\section{Introduction} \label{section1}
The path to understanding the internal structure of geometries is through understanding the elementary notions
and elementary geometry of the space under study. Thus, mapping the structure of Thurston geometries can also provide
important information. One of the most important results may be the knowing of non-Euclidean crystallography, since
it cannot be ruled out that under certain extreme conditions, materials
crystallize in a ``different geometric structure´´.

In the study of the internal structures of Thurston spaces, it is worth distinguishing 
two significantly different directions related to the translation distance or geodesic distance.

We can introduce in a natural way (see \cite{M97})
translation mappings any point to any other point. Consider a unit tangent
vector at the origin. Translations carry this vector to a tangent vector any
other point.
If a curve $t\rightarrow (x(t),y(t),z(t))$ has just the translated vector as its tangent vector at
each point, then the curve is called a {\it translation curve}. This assumption leads to
a system of first order differential equations. Thus translation curves are simpler
than geodesics and differ from them in $\NIL$, $\SLR$ and $\SOL$ geometries.
In $\EUC$, $\SPH$, $\HYP$, $\SXR$ and $\HXR$ geometries, the translation and geodesic curves coincide
with each other. But in the $\NIL$, $\SLR$ and $\SOL$ geometries,
translation curves are in many ways more natural than geodesics. Therefore, we
distinguish two different distance functions: $d^g$ is the usual geodesic distance function,
and $d^t$ is the translation distance function. So we obtain
two types of curves, triangles, bisector surfaces (and two types of the corresponding Dirichlet-Voronoi cells) etc. from the two different distance functions,
but {\it{in the present paper we consider only the translation case}}. 

In our previous works we have dealt with both cases and since it plays 
a fundamental role in the structure of geometries, 
we first dealt with the concept and properties of the translation and geodesic triangles.
The extension of these notions are interesting in Thurston geometries 
with non-constant curvature, since the structures of Euclidean $\EUC$, spherical $\SPH$ and hyperbolic $\HYP$ spaces are well known. 

In each Thurston space we have examined the possible values of the internal 
angle sums of the geodesic and translational triangles (see \cite{CsSz16}, \cite{Cs-Sz24}, \cite{Sz18}, \cite{Sz20}, \cite{Sz22-1}).

A second important direction is the study of Apollonius surfaces, a special case of 
which is the bisector surfaces. These are of fundamental importance due to the determinations of Dirichlet-Voronoi cells.
The geodesic-like Apollonius surfaces are investigated in $\SXR$, $\HXR$ and $\NIL$ geometries in \cite{Sz22-2, Sz23-1} and in \cite{Cs-Sz25} 
the translation-like Apollonius surfaces in Thurston geometries.
Moreover, we studied the translation-like equidistant surfaces in $\SOL$ and $\NIL$ geometries in \cite{Sz19} and in \cite{VSz19}.
In \cite{PSSz10}, \cite{PSSz11-1}, \cite{PSSz11-2} we studied the geodesic-like bisector surfaces in $\SXR$ and $\HXR$ spaces.
\begin{rmrk}
A special case of Apollonius surfaces is the geodesic- or translation-like bisector (or equidistant) surface 
of two arbitrary points. These surfaces have an important role in structure of Dirichlet-Voronoi (briefly, D-V) cells.
The D-V-cells are relevant in the study of tilings, ball packing and ball covering. E.g. if the point set is the orbit of a point - generated by
a discrete isometry group of the considered space - then we obtain a monohedral D-V cell decomposition (tiling) of the considered space and it is interesting to examine its
optimal ball packing and covering. In $3$-dimensional spaces of constant curvature, the
D-V cells have been widely investigated, but in the other Thurston geometries $\SXR$, $\HXR$,
$\NIL$, $\SOL$, $\SLR$ there are few results on this topic.
\end{rmrk}
We examined further interesting problems related to the Thurston geometries 
with non-constant curvature in papers \cite{MSz12, MSzV17, Sz07, Sz10, Sz11, Sz12, 
Sz14-1, Sz14-2, Sz23-2, Sz24}.

The next important question, which is also closely related to the
Apollonius surfaces, is the determination of the surface of the given geodesic or translational triangle. 
Defining this is an essential condition for stating elementary geometric concepts and theorems related to triangles.

However, defining the surface of a translation or geodesic triangle in Thurston geometries with non-constant curvature is not 
straightforward. The usual geodesic or translation-like triangle surface definition is not possible because 
the geodesic or translation curves starting from different vertices and ending at points of the corresponding opposite edges define different
surfaces in general, i.e. {\it geodesics or translation curves starting from different vertices and ending at points on the corresponding opposite side usually do not intersect.}

In the works, \cite{Sz22-2}, \cite{Sz23-1} we proposed a possible definition of geodesic-like triangular surfaces in $\SXR$, $\HXR$ and $\NIL$ geometries, 
which was done with the help of Apollonius surfaces and returned traditional triangular surfaces in geometries with constant curvature.

In the paper \cite{Cs-Sz25}, we extended these questions to the translation triangles of non-constant curvature Thurston geometries. 
We provided a new possible definition of the surfaces of translation-like triangles and determined their equations.
In the following we will use these definitions of the surfaces of translation triangles in Thurston geometrien.

In Thurston geometries of constant curvature the Ceva's and Menelaus' theorems are well known and in \cite{Sz22-2} and  \cite{Sz23-1} we generalized 
the Ceva's and Menelaus' theorems to $\SXR$, $\HXR$ and $\NIL$ spaces.

In this paper we will continue to investigate this issue, generalize and prove the theorems of Menelaus and Ceva for translation triangles in $\NIL$, $\SLR$ and $\SOL$ spaces. 

The computations and the proofs are based on the projective models 
of Thurston geometries described by E. Moln\'ar in \cite{M97}.
\section{Projective model of $\NIL$, $\SLR$ and $\SOL$ geometries}
In the following, we summarize the most important concepts of the further investigated geometries.
\subsection{On $\NIL$ geometry}
$\NIL$ geometry can be derived from the famous real matrix group $\mathbf{L(\mathbf{R})}$ discovered by Werner Heisenberg. The left (row-column) 
multiplication of Heisenberg matrices
     \begin{equation}
     \begin{gathered}
     \begin{pmatrix}
         1&x&z \\
         0&1&y \\
         0&0&1 \\
       \end{pmatrix}
       \begin{pmatrix}
         1&a&c \\
         0&1&b \\
         0&0&1 \\
       \end{pmatrix}
       =\begin{pmatrix}
         1&a+x&c+xb+z \\
         0&1&b+y \\
         0&0&1 \\
       \end{pmatrix}
      \end{gathered} \label{2.7}
     \end{equation}
defines "translations" $\mathbf{L}(\mathbf{R})= \{(x,y,z): x,~y,~z\in \mathbf{R} \}$ 
on the points of $\NIL= \{(a,b,c):a,~b,~c \in \mathbf{R}\}$. 
These translations are not commutative in general. The matrices $\mathbf{K}(z) \vartriangleleft \mathbf{L}$ of the form
     \begin{equation}
     \begin{gathered}
       \mathbf{K}(z) \ni
       \begin{pmatrix}
         1&0&z \\
         0&1&0 \\
         0&0&1 \\
       \end{pmatrix}
       \mapsto (0,0,z)  
      \end{gathered}\label{2.8}
     \end{equation} 
constitute the one parametric centre, i.e. each of its elements commutes with all elements of $\mathbf{L}$. 
The elements of $\mathbf{K}$ are called {\it fibre translations}. $\NIL$ geometry of the Heisenberg group can be projectively 
(affinely) interpreted by "right translations" 
on points as the matrix formula 
     \begin{equation}
     \begin{gathered}
       (1;a,b,c) \to (1;a,b,c)
       \begin{pmatrix}
         1&x&y&z \\
         0&1&0&0 \\
         0&0&1&x \\
         0&0&0&1 \\
       \end{pmatrix}
       =(1;x+a,y+b,z+bx+c) 
      \end{gathered} \label{2.9}
     \end{equation} 
shows, according to (\ref{2.7}). Here we consider $\mathbf{L}$ as projective collineation group with right actions in homogeneous coordinates.
We will use the usual projective model of $\NIL$ (see \cite{M97} and \cite{MSz06}).

The translation group $\mathbf{L}$ defined by formula (\ref{2.9}) can be extended to a larger group $\mathbf{G}$ of collineations,
preserving the fibres, that will be equivalent to the (orientation preserving) isometry group of $\NIL$. 

In \cite{MSz06} we has shown that 
a rotation through angle $\omega$
about the $z$-axis at the origin, as isometry of $\NIL$, keeping invariant the Riemann
metric everywhere, will be a quadratic mapping in $x,y$ to $z$-image $\overline{z}$ as follows:
     \begin{equation}
     \begin{gathered}
       \mathcal{M}=\br(O,\omega):(1;x,y,z) \to (1;\overline{x},\overline{y},\overline{z}); \\ 
       \overline{x}=x\cos{\omega}-y\sin{\omega}, \ \ \overline{y}=x\sin{\omega}+y\cos{\omega}, \\
       \overline{z}=z-\frac{1}{2}xy+\frac{1}{4}(x^2-y^2)\sin{2\omega}+\frac{1}{2}xy\cos{2\omega}.
      \end{gathered} \label{2.10}
     \end{equation}
This rotation formula $\mathcal{M}$, however, is conjugate by the quadratic mapping $\alpha$ to the linear rotation $\Omega$ as follows
     \begin{equation}
     \begin{gathered}
       \alpha^{-1}: \ \ (1;x,y,z) \stackrel{\alpha^{-1}}{\longrightarrow} (1; x',y',z')=(1;x,y,z-\frac{1}{2}xy) \ \ \text{to} \\
       \Omega: \ \ (1;x',y',z') \stackrel{\Omega}{\longrightarrow} (1;x",y",z")=(1;x',y',z')
       \begin{pmatrix}
         1&0&0&0 \\
         0&\cos{\omega}&\sin{\omega}&0 \\
         0&-\sin{\omega}&\cos{\omega}&0 \\
         0&0&0&1 \\
       \end{pmatrix}, \\
       \text{with} \ \ \alpha: (1;x",y",z") \stackrel{\alpha}{\longrightarrow}  (1; \overline{x}, \overline{y},\overline{z})=(1; x",y",z"+\frac{1}{2}x"y").
      \end{gathered} \label{2.11}
     \end{equation}
This quadratic conjugacy modifies the $\NIL$ translations in (2.3), as well. 
\subsubsection{Translation curves}
We consider a $\NIL$ curve $(1,x(t), y(t), z(t) )$ with a given starting tangent vector at the origin $O=E_0=(1,0,0,0)$
\begin{equation}
   \begin{gathered}
      u=\dot{x}(0),\ v=\dot{y}(0), \ w=\dot{z}(0).
       \end{gathered} \label{2.12}
     \end{equation}
For a translation curve let its tangent vector at the point $(1,x(t), y(t), z(t) )$ be defined by the matrix (\ref{2.9}) 
with the following equation:
\begin{equation}
     \begin{gathered}
     (0,u,v,w)
     \begin{pmatrix}
         1&x(t)&y(t)&z(t) \\
         0&1&0&0 \\
         0&0&1&x(t) \\
         0&0&0&1 \\
       \end{pmatrix}
       =(0,\dot{x}(t),\dot{y}(t),\dot{z}(t)).
       \end{gathered} \label{2.13}
     \end{equation}
Thus, the {\it translation curves} in $\NIL$ geometry (see  \cite{Cs-Sz23}, \cite{MSz06} \cite{MSzV}) are defined by the above first order differential equation system 
$\dot{x}(t)=u, \ \dot{y}(t)=v,  \ \dot{z}(t)=v \cdot x(t)+w,$ whose solution is the following: 
\begin{equation}
   \begin{gathered}
       x(t)=u t, \ y(t)=v t,  \ z(t)=\frac{1}{2}uvt^2+wt.
       \end{gathered} \label{2.14}
\end{equation}
We assume that the starting point of a translation curve is the origin, because we can transform a curve into an 
arbitrary starting point by translation (2.3), moreover, unit initial velocity translation can be assumed by "geographic" parameters $\phi$ and $\theta$:
\begin{equation}
\begin{gathered}
        x(0)=y(0)=z(0)=0; \\ \ u=\dot{x}(0)=\cos{\theta} \cos{\phi}, \ \ v=\dot{y}(0)=\cos{\theta} \sin{\phi}, \ \ w=\dot{z}(0)=\sin{\theta}; \\ 
        - \pi \leq \phi \leq \pi, \ -\frac{\pi}{2} \leq \theta \leq \frac{\pi}{2}. \label{2.15}
\end{gathered}
\end{equation}
\subsection{On Sol geometry}
\label{sec:1}
In this Section we summarize the significant notions and notations of real $\SOL$ geometry (see \cite{M97}, \cite{S}).

$\SOL$ is defined as a 3-dimensional real Lie group with multiplication
\begin{equation}
     \begin{gathered}
(a,b,c)(x,y,z)=(x + a e^{-z},y + b e^z ,z + c).
     \end{gathered} \label{2.16}
     \end{equation}
We note that the conjugacy by $(x,y,z)$ leaves invariant the plane $(a,b,c)$ with fixed $c$:
\begin{equation}
     \begin{gathered}
(x,y,z)^{-1}(a,b,c)(x,y,z)=(x(1-e^{-c})+a e^{-z},y(1-e^c)+b e^z ,c).
     \end{gathered} \label{2.17}
     \end{equation}
Moreover, for $c=0$, the action of $(x,y,z)$ is only by its $z$-component, where $(x,y,z)^{-1}=(-x e^{z}, -y e^{-z} ,-z)$. Thus the $(a,b,0)$ plane is distinguished as a {\it base plane} in
$\SOL$, or by other words, $(x,y,0)$ is normal subgroup of $\SOL$.
$\SOL$ multiplication can also be affinely (projectively) interpreted by ``right translations"
on its points as the following matrix formula shows, according to (\ref{2.16}):
     \begin{equation}
     \begin{gathered}
     (1,a,b,c) \to (1,a,b,c)
     \begin{pmatrix}
         1&x&y&z \\
         0&e^{-z}&0&0 \\
         0&0&e^z&0 \\
         0&0&0&1 \\
       \end{pmatrix}
       =(1,x + a e^{-z},y + b e^z ,z + c)
       \end{gathered} \label{2.18}
     \end{equation}
by row-column multiplication.
This defines ``translations" $\mathbf{L}(\mathbf{R})= \{(x,y,z): x,~y,~z\in \mathbf{R} \}$
on the points of space $\SOL= \{(a,b,c):a,~b,~c \in \mathbf{R}\}$.
These translations are not commutative, in general.
Here we can consider $\mathbf{L}$ as projective collineation group with right actions in homogeneous
coordinates as usual in classical affine-projective geometry.

We will use the usual projective model of $\SOL$ (see \cite{M97} and \cite{Sz19}).

It will be important for us that the full isometry group Isom$(\SOL)$ has eight components, since the stabilizer of the origin
is isomorphic to the dihedral group $\mathbf{D_4}$, generated by two involutive (involutory) transformations:
\begin{equation}
   \begin{gathered}
      (1)  \ \ y \leftrightarrow -y; \ \ (2)  \ x \leftrightarrow y; \ \ z \leftrightarrow -z; \ \ \text{i.e. first by $3\times 3$ matrices}:\\
     (1) \ \begin{pmatrix}
               1&0&0 \\
               0&-1&0 \\
               0&0&1 \\
     \end{pmatrix}; \ \ \
     (2) \ \begin{pmatrix}
               0&1&0 \\
               1&0&0 \\
               0&0&-1 \\
     \end{pmatrix}; \\
     \end{gathered} \label{2.19}
     \end{equation}
     with its product, generating a cyclic group $\mathbf{C_4}$ of order 4
     \begin{equation}
     \begin{gathered}
     \begin{pmatrix}
                    0&1&0 \\
                    -1&0&0 \\
                    0&0&-1 \\
     \end{pmatrix};\ \
     \begin{pmatrix}
               -1&0&0 \\
               0&-1&0 \\
               0&0&1 \\
     \end{pmatrix}; \ \
     \begin{pmatrix}
               0&-1&0 \\
               1&0&0 \\
               0&0&-1 \\
     \end{pmatrix};\ \
     \mathbf{Id}=\begin{pmatrix}
               1&0&0 \\
               0&1&0 \\
               0&0&1 \\
     \end{pmatrix}.
     \end{gathered} \notag
     \end{equation}
     Or we write by collineations fixing the origin $O=(1,0,0,0)$:
\begin{equation}
(1) \ \begin{pmatrix}
         1&0&0&0 \\
         0&1&0&0 \\
         0&0&-1&0 \\
         0&0&0&1 \\
       \end{pmatrix}, \ \
(2) \ \begin{pmatrix}
         1&0&0&0 \\
         0&0&1&0 \\
         0&1&0&0 \\
         0&0&0&-1 \\
       \end{pmatrix} \ \ \text{of form (2.13)}. \label{2.20}
\end{equation}
A general isometry of $\SOL$ to the origin $O$ is defined 
by a product $\gamma_O \tau_X$, first $\gamma_O$ of form (2.14) then 
$\tau_X$ of (2.12). To
a general point $A=(1,a,b,c)$, this will be a product $\tau_A^{-1} \gamma_O \tau_X$, mapping $A$ into $X=(1,x,y,z)$.

We remark only that the role of $x$ and $y$ can be exchanged throughout the paper, but this leads to the mirror interpretation of $\SOL$.
As formula (2.10) fixes the metric of $\SOL$, the change above is not an isometry of a fixed $\SOL$ interpretation. Other conventions are also accepted
and used in the literature.

{\it $\SOL$ is an affine metric space (affine-projective one in the sense of the unified formulation of \cite{M97}). Therefore, its linear, affine, unimodular,
etc. transformations are defined as those of the embedding affine space.}
\subsubsection{Translation curves}
We consider a $\SOL$ curve $(1,x(t), y(t), z(t) )$ with a given starting tangent vector at the origin $O=(1,0,0,0)$
\begin{equation}
   \begin{gathered}
      u=\dot{x}(0),\ v=\dot{y}(0), \ w=\dot{z}(0).
       \end{gathered} \label{2.21}
     \end{equation}
For a translation curve let its tangent vector at the point $(1,x(t), y(t), z(t) )$ be defined by the matrix (\ref{2.18})
with the following equation:
\begin{equation}
     \begin{gathered}
     (0,u,v,w)
     \begin{pmatrix}
         1&x(t)&y(t)&z(t) \\
         0&e^{-z(t)}&0&0 \\
         0&0&e^{z(t)}& 0 \\
         0&0&0&1 \\
       \end{pmatrix}
       =(0,\dot{x}(t),\dot{y}(t),\dot{z}(t)).
       \end{gathered} \label{2.22}
     \end{equation}
Thus, {\it translation curves} in $\SOL$ geometry (see \cite{Sz19} and \cite{Sz22-1}) are defined by the first order differential equation system
$\dot{x}(t)=u e^{-z(t)}, \ \dot{y}(t)=v e^{z(t)},  \ \dot{z}(t)=w,$ whose solution is the following:
\begin{equation}
   \begin{gathered}
     x(t)=-\frac{u}{w} (e^{-wt}-1), \ y(t)=\frac{v}{w} (e^{wt}-1),  \ z(t)=wt, \ \mathrm{if} \ w \ne 0 \ \mathrm{and} \\
     x(t)=u t, \ y(t)=v t,  \ z(t)=z(0)=0 \ \ \mathrm{if} \ w =0.
       \end{gathered} \label{2.23}
\end{equation}
We assume that the starting point of a translation curve is the origin, because we can transform a curve into an
arbitrary starting point by translation (\ref{2.18}), moreover, unit velocity translation can be assumed :
\begin{equation}
\begin{gathered}
        x(0)=y(0)=z(0)=0; \\ \ u=\dot{x}(0)=\cos{\theta} \cos{\phi}, \ \ v=\dot{y}(0)=\cos{\theta} \sin{\phi}, \ \ w=\dot{z}(0)=\sin{\theta}; \\
        - \pi < \phi \leq \pi, \ -\frac{\pi}{2} \leq \theta \leq \frac{\pi}{2}. \label{2.24}
\end{gathered}
\end{equation}
Thus we obtain the parametric equation of the {\it translation curve segment} $t(\phi,\theta,t)$ with starting point at the origin in direction
\begin{equation}
\bt(\phi, \theta)=(\cos{\theta} \cos{\phi}, \cos{\theta} \sin{\phi}, \sin{\theta}) \label{2.25}
\end{equation}
where $t \in [0,r] ~ r \in \bR^+$. If $\theta \ne 0$ then the system of equation is:
\begin{equation}
\begin{gathered}
        \left\{ \begin{array}{ll}
        x(\phi,\theta,t)=-\cot{\theta} \cos{\phi} (e^{-t \sin{\theta}}-1), \\
        y(\phi,\theta,t)=\cot{\theta} \sin{\phi} (e^{t \sin{\theta}}-1), \\
        z(\phi,\theta,t)=t \sin{\theta}.
        \end{array} \right. \\
        \text{If $\theta=0$ then}: ~  x(t)=t\cos{\phi} , \ y(t)=t \sin{\phi},  \ z(t)=0.
        \label{2.26}
        \end{gathered}
\end{equation}
\subsection{On $\SLR$ geometry}
The real $ 2\times 2$ matrices $\begin{pmatrix}
         d&b \\
         c&a \\
         \end{pmatrix}$ with unit determinant $ad-bc=1$
constitute a Lie transformation group by the usual product operation, taken to act on row matrices as on point coordinates on the right as follows
\begin{equation}
\begin{gathered}
(z^0,z^1)\begin{pmatrix}
         d&b \\
         c&a \\
         \end{pmatrix}=(z^0d+z^1c, z^0 b+z^1a)=(w^0,w^1)\\
\mathrm{with} \ w=\frac{w^1}{w^0}=\frac{b+\frac{z^1}{z^0}a}{d+\frac{z^1}{z^0}c}=\frac{b+za}{d+zc}, \label{2.27}
\end{gathered}
\end{equation}
as action on the complex projective line $\bC^{\infty}$ (see \cite{M97}, \cite{MSz}).
This group is a $3$-dimensional manifold, because of its $3$ independent real coordinates and with its usual neighbourhood topology (\cite{S}, \cite{T}).
In order to model the above structure in the projective sphere $\cP \cS^3$ and in the projective space $\cP^3$ (see \cite{M97}),
we introduce the new projective coordinates $(x^0,x^1,x^2,x^3)$ where
$a:=x^0+x^3, \ b:=x^1+x^2, \ c:=-x^1+x^2, \ d:=x^0-x^3$
with the positive, then the non-zero multiplicative equivalence as projective freedom in $\cP \cS^3$ and in $\cP^3$, respectively.
Then it follows that $0>bc-ad=-x^0x^0-x^1x^1+x^2x^2+x^3x^3$
describes the interior of the above one-sheeted hyperboloid solid $\cH$ in the usual Euclidean coordinate simplex with the origin
$E_0=(1,0,0,0)$ and the ideal points of the axes $E_1^\infty=(0,1,0,0)$, $E_2^\infty(0,0,1,0)$, $E_3^\infty=(0,0,0,1)$.
We consider the collineation group ${\bf G}_*$ that acts on the projective sphere $\cS\cP^3$  and preserves a polarity i.e. a scalar product of signature
$(- - + +)$, this group leaves the one sheeted hyperboloid solid $\cH$ invariant.
We have to choose an appropriate subgroup $\mathbf{G}$ of $\mathbf{G}_*$ as isometry group, then the universal covering group and space
$\widetilde{\cH}$ of $\cH$ will be the hyperboloid model of $\SLR$ \cite{M97}.

The elements of the isometry group of
$\mathbf{SL_2R}$ (and so by the above extension the isometries of $\SLR$) can be described in \cite{M97} and \cite{MSz}).
Moreover, we have the projective proportionality, of course.
We define the {\it translation group} $\bG_T$, as a subgroup of the isometry group of $\mathbf{SL_2R}$,
the isometries acting transitively on the points of ${\cH}$ and by the above extension on the points of $\SLR$ and $\widetilde{\cH}$.
$\bG_T$ maps the origin $E_0=(1,0,0,0)$ onto $X=(x^0,x^1,x^2,x^3)$. These isometries and their inverses (up to a positive determinant factor)
are given by the following matrices:
\begin{equation}
\begin{gathered} \bT:~(t_i^j)=
\begin{pmatrix}
x^0&x^1&x^2&x^3 \\
-x^1&x^0&x^3&-x^2 \\
x^2&x^3&x^0&x^1 \\
x^3&-x^2&-x^1&x^0
\end{pmatrix}.
\end{gathered} \label{2.28}
\end{equation}
The rotation about the fibre line through the origin $E_0=(1,0,0,0)$ by angle $\omega$ $(-\pi<\omega\le \pi)$ can be expressed by the following matrix
(see \cite{M97})
\begin{equation}
\begin{gathered} \bR_{E_O}(\omega):~(r_i^j(E_0,\omega))=
\begin{pmatrix}
1&0&0&0 \\
0&1&0&0 \\
0&0&\cos{\omega}&\sin{\omega} \\
0&0&-\sin{\omega}&\cos{\omega}
\end{pmatrix},
\end{gathered} \label{2.29}
\end{equation}
and the rotation $\bR_X(\omega)$ about the fibre line through $X=(x^0,x^1,x^2,x^3)$ by angle $\omega$ can be derived by formulas (\ref{2.28}) and (\ref{2.29}):
\begin{equation}
\bR_X(\omega)=\bT^{-1} \bR_{E_O} (\omega) \bT:~(r_i^j(X,\omega)).
\label{2.30}
\end{equation}
After \cite{M97}, we introduce the so-called hyperboloid parametrization as follows
\begin{equation}
\begin{gathered}
x^0=\cosh{r} \cos{\phi}, ~ ~
x^1=\cosh{r} \sin{\phi}, \\
x^2=\sinh{r} \cos{(\theta-\phi)}, ~ ~
x^3=\sinh{r} \sin{(\theta-\phi)},
\end{gathered} \label{2.31}
\end{equation}
where $(r,\theta)$ are the polar coordinates of the base plane and $\phi$ is just the fibre coordinate. We note that
$$-x^0x^0-x^1x^1+x^2x^2+x^3x^3=-\cosh^2{r}+\sinh^2{r}=-1<0.$$
The inhomogeneous coordinates corresponding to (2.25), that play an important role in the later visualization of prism tilings in $\EUC$,
are given by
\begin{equation}
\begin{gathered}
x=\frac{x^1}{x^0}=\tan{\phi}, ~ ~
y=\frac{x^2}{x^0}=\tanh{r} \frac{\cos{(\theta-\phi)}}{\cos{\phi}}, \\
z=\frac{x^3}{x^0}=\tanh{r} \frac{\sin{(\theta-\phi)}}{\cos{\phi}}.
\end{gathered} \label{2.32}
\end{equation}
\subsubsection{Translation curves}
We recall some basic facts about translation curves in $\SLR$  following 
\cite{MSzV, MSz, MSzV17}. For any point $X=(x^{0}, x^{1}, x^{2}, x^{3}) \in {\mathcal H}$ 
(and later also for points in $\widetilde{\mathcal H}$) the \emph{translation map} from  the origin 
$E_{0}=(1;0;0;0)$ to $X$ is defined by the \emph{translation matrix} $\bT$ and its inverse presented in (\ref{2.28}).

Let us consider for a given vector  $(q;u;v;w)$ a curve 
$\mathcal C(t) = (x^{0}(t,;x^{1}(t),$ $x^{2}(t), x^{3}(t))$, $t \geq 0$, in ${\mathcal H}$ starting at the origin: $\mathcal C(0) = E_{0}(1;0;0;0)$ and such that
$$
\dot{\mathcal C} (0) = (\dot{x}^{0}(0), \dot{x}^{1}(0), \dot{x}^{2}(0), \dot{x}^{3}(0)) = (q,u,v,w),
$$
where $\dot{\mathcal C}(t)  = (\dot{x}^{0}(t), \dot{x}^{1}(t), \dot{x}^{2}(t), \dot{x}^{3}(t))$
is the tangent vector at any point of the curve. For $t \geq 0$ there exists a matrix
\begin{equation*}
\bT(t) = \left( \begin{array}{cccc}
x^{0}(t) & x^{1}(t) & x^{2}(t) & x^{3}(t) \cr
-x^{1}(t) & x^{0}(t) & x^{3}(t) & -x^{2}(t) \cr
x^{2}(t) & x^{3}(t) & x^{0}(t) & x^{1}(t) \cr
x^{3}(t) & -x^{2}(t) & -x^{1}(t) & x^{0}(t) \cr
\end{array} \right) \tag{2.27}
\end{equation*}
which defines the translation from $\mathcal C(0)$ to $\mathcal C (t)$:
\begin{equation*}
\mathcal C(0) \cdot \bT(t) = \mathcal C(t),  \quad t \geqslant 0. \tag{2.28}
\end{equation*}
The $t$-parametrized family $\bT(t)$ of translations is used in the following definition.

As we mentioned this earlier, 
the curve $\mathcal C(t)$, $t \geqslant 0$, is said to be a \emph{translation curve} if
\begin{equation*}
\dot{\mathcal C }(0) \cdot \bT(t) = \dot{\mathcal C}(t),  \quad t \geqslant 0. \tag{2.29}
\end{equation*}

The solution, depending on $(q,u,v,w)$ can be seen in \cite{Sz19}, where it splits into three cases. 

It was observed above that for any $X=(x^{0},x^{1},x^{2},x^{3}) \in \widetilde{\mathcal H}$ there 
is a suitable transformation $\bT^{-1}$, which sent $X$ to the origin $E_{0}$ along a translation curve.
For a given translation curve $\mathcal C = \mathcal C (t)$ the initial unit tangent vector $(u,v,w)$ (in Euclidean coordinates) 
at $E_{0}$ can be presented as
\begin{equation}
u = \sin \alpha, \quad v = \cos \alpha \cos \lambda, \quad w = \cos \alpha \sin \lambda, \notag
\end{equation}
for some $-\frac{\pi}{2} \leqslant \alpha \leqslant \frac{\pi}{2}$ and $ -\pi < \lambda \leqslant \pi$.  
In $\widetilde{\mathcal H}$ this vector is of length square $-u^{2} + v^{2} + w^{2}  = \cos 2 \alpha$. We 
always can assume that $\mathcal C$ is parametrized by the translation arc-length parameter $t = s \geqslant 0$. 
Then coordinates of a point $X=(x; y; z)$ of $\mathcal C$, such that the translation distance between $E_{0}$ and $X$ equals 
$s$, depend on $(\lambda,\alpha,s)$ as geographic coordinates according to the above considered three cases as follows.
\begin{table}[ht]
\caption{Translation  curves.} \label{table2}
\vspace{3mm}
\centerline{
$
\begin{array}{|c|l|} \hline
\textrm{direction} & \textrm{parametrization of a translation curve}  \cr \hline  
\begin{gathered}
0 \le \alpha < \frac{\pi}{4} \\ (\bH^2-{\rm like})
\end{gathered}
&
\begin{array}{l}
\begin{gathered} \noalign{\vskip2pt}\begin{pmatrix} x(s,\alpha,\lambda)\\y(s,\alpha,\lambda)\\z(s,\alpha,\lambda)\end{pmatrix}=  \frac{\tanh (s \sqrt{\cos 2 \alpha})}{\sqrt{\cos 2\alpha}} \begin{pmatrix} \sin \alpha\\
 \cos \alpha \cos \lambda\\
 \cos \alpha \sin \lambda)\end{pmatrix} \end{gathered} \cr \noalign{\vskip2pt}\end{array} \cr
 \hline
\begin{gathered} \alpha=\frac{\pi}{4} \\ ({\rm light-like}) \end{gathered}
&
\begin{array}{l}
\begin{gathered} \noalign{\vskip2pt}\begin{pmatrix} x(s,\alpha,\lambda)\\y(s,\alpha,\lambda)\\z(s,\alpha,\lambda)\end{pmatrix}=  \frac{\sqrt{2} s}{2} 
\begin{pmatrix} 1\\
 \cos \lambda\\
 \sin \lambda)\end{pmatrix} \end{gathered} \cr
\noalign{\vskip2pt}\end{array} \cr \hline
\begin{gathered} \frac{\pi}{4}  < \alpha \le \frac{\pi}{2} \\ ({\rm fibre-like}) \end{gathered}
&
\begin{array}{l}  \begin{gathered} \noalign{\vskip2pt}\begin{pmatrix} x(s,\alpha,\lambda)\\y(s,\alpha,\lambda)\\z(s,\alpha,\lambda)\end{pmatrix}=  \frac{\tan (s \sqrt{-\cos 2 \alpha})}{\sqrt{-\cos 2\alpha}} \begin{pmatrix} \sin \alpha\\
 \cos \alpha \cos \lambda\\
 \cos \alpha \sin \lambda)\end{pmatrix} \end{gathered}   \cr
\noalign{\vskip2pt}
\end{array}  \cr \hline
\end{array}
$
}
\end{table}
\section{Menelaus' and Ceva's theorems in $\NIL$, $\SOL$ and $\SLR$ spaces.}
\subsection{The surfaces of translation triangles}
We consider $3$ points $A_0$, $A_1$, $A_2$ in the projective model of the space $X$ (see Section 2) $(X\in\{\NIL, \SOL, \SLR \})$.
The {\it translation segments} $a_k$ connecting the points $A_i$ and $A_j$
$(i<j,~i,j,k \in \{0,1,2\}, k \ne i,j$) are called sides of the {\it translation triangle} with vertices $A_0$, $A_1$, $A_2$.

However, defining the surface of a translation triangle in the spaces $\NIL$ and $\SOL$ is not straightforward. The usual translation-like
triangle surface definition in these geometries is
not possible because the translation curves starting from different vertices and ending at points of the corresponding opposite edges define different
surfaces, i.e. {\it translation curves starting from different vertices and ending at points on the corresponding opposite side usually do not intersect.} 
In $\SLR$ geometry, this problem does not arise because translation curves are lines in the Euclidean sense (see Table 1 in Section 2), 
so the surface of the translation triangle can be defined by the Euclidean plane given by its side lines.

Therefore, we introduced a new definition of the surface $\mathcal{S}_{A_0A_1A_2}$ of the translation triangle in \cite{Cs-Sz25}: 
{\it The $S^{X,t}_{A_0A_1A_2}$ translation-like triangular surface of translation triangle $A_0A_1A_2$ in the Thurston geometry $X \in \{\EUC,\SPH,\HYP,\SXR,\HXR,\NIL,\SLR,\SOL\}$ is  
the set of all $P$ points of $X$ from which the tangents in $P$ of the translation curves drawn to vertices $A_0$ $A_1$ and $A_2$ are coplanar.}
We also determined the equations of these surfaces and visualized them and we will use this surface definition in the following.
An important question is how to define the curves that pass through the surface of a translation triangle. 
Their definitions always depend on the structure of the space, so their definition is always done when examining the given space.

\subsection{Known results in $\EUC$, $\SPH$ $\HYP$, $\SXR$ and $\HXR$.}
The Ceva's and Menelaus' theorems in Euclidean space are well known, here we only mention the results for the additional spaces.
\begin{enumerate}
\item {\bf Sperical and hyperbolic cases}

First we recall the definition of simple ratios in the sphere $\bS^2$ and the
plane $\bH^2$.
The models of the above plane geometries of constant curvature are
embedded in the models of the previously described geometries $\SXR$ and $\HXR$
as ``base planes" and are used hereinafter for our discussions.

A spherical triangle is the space included by arcs of great circles on the surface of a sphere, subject to the limitation that these arcs and the
further circle arcs in the following items in the spherical
plane are always less or equal than a semicircle.
\begin{definition}
If $A$, $B$ and $P$ are distinct points on a line in
the $Y\in\{\EUC,\HYP,\SPH\}$ space, then
their simple ratio is
$$s_g^Y(A,P,B) =  w^Y(d^Y(A,P))/w^Y(d^Y(P,B)),$$ if $P$ is between $A$ and $B$, and
$$s_g^Y(A,P,B) = -w^Y(d^Y(A,P))/w^Y(d^Y(P,B)),$$ otherwise where
$w^Y(x):=x$ if $Y=\EUC$, $w^Y(x):=sin(x)$ if $Y=\SPH$ and $w^Y(x):=sinh(x)$ if $Y=\HYP$.
\end{definition}

With this definition, the corresponding sine rule of the geometry $Y$ leads to
Menelaus's and Ceva's theorems \cite{K,PS14}:

\begin{Theorem} [Menelaus's Theorem for triangles in the space $Y$]
If is a line $l$ not through any vertex of an triangle $ABC$ such that
$l$ meets $BC$ in $Q$, $AC$ in $R$, and $AB$ in $P$,
then $$s^Y(A,P,B)s^Y(B,Q,C)s^Y(C,R,A) = -1.$$ ~ ~ $\square$
\end{Theorem}

\begin{Theorem}[Ceva's Theorem for triangles in the space $Y$]
If $T$ is a point not on any side of a triangle $ABC$ such that
$AT$ and $BC$ meet in $Q$, $BX$ and $AC$ in $R$, and $CX$ and $AB$ in $P$,
then $$s^Y(A,P,B)s^Y(B,Q,C)s^Y(C,R,A) = 1.$$ ~ ~ $\square$
\end{Theorem}
\begin{rmrk}
It is easy to see that the ``reversals'' of the above theorems are also true.
\end{rmrk}
\item {\bf $\SXR$ and $\HXR$ cases} 

We extended in \cite{Sz22-2} the definition of the simple ratio to the $X\in\{\SXR, \HXR \}$ spaces.
If $X=\SXR$ then is clear that the space contains its ``base sphere" which is a geodesic surface (or translation-like, here they coincide).
Therefore, similarly to the spherical spaces we assume that the geodesic arcs
in the following items are always less than or equal to a semicircle.
\begin{definition}
If $A$, $B$ and $P$ are distinct points on a non-fibrum-like geodesic curve in the $X \in \{\SXR, \HXR \}$ space , then
their simple ratio is
$$s_g^X(A,P,B) =  w^X\Big({d^X(A,P)}{\cos(v)}\Big)/w^X\Big({d^X(P,B)}{\cos(v)}\Big),$$ if $P$ is between $A$ and $B$, and
$$s_g^X(A,P,B) = -w^X\Big({d^X(A,P)}{\cos(v)}\Big)/w^X\Big({d^X(P,B)}{\cos(v)}\Big),$$
otherwise where
$w^X(x):=sin(x)$ if $X=\SXR$, $w^Y(x):=sinh(x)$ if $X=\HXR$ and $v$ is
the parameter of the geodesic curve containing points $A,B,P$.
\end{definition}
\begin{Theorem}[Ceva's Theorem for triangles in general location, \cite{Sz22-2}]
If $T$ is a point not on any side of a geodesic triangle $A_0A_1A_2$ in $X\in\{\SXR, \HXR\}$ such that
the curves $A_0T$ and $g_{A_1A_2}^X$ meet in $Q$, $A_1T$ and $g_{A_0A_2}^X$ in $R$, and $A_2T$ and $g_{A_0A_1}^X$ in $P$, $(A_0T, A_1T, A_2T \subset \mathcal{S}^X_{A_0A_1A_2})$
then $$s_g^X(A_0,P,A_1)s_g^X(A_1,Q,A_2)s_g^X(A_2,R,A_0) = 1.$$
\end{Theorem}
\begin{Theorem}[Menelaus's theorem for triangles in general location, \cite{Sz22-2}]
If $l$ is a line not through any vertex of a geodesic triangle
$A_0A_1A_2$ lying in a surface $\mathcal{S}^X_{A_0A_1A_2}$ in the $X\in\{\SXR, \HXR\}$
geometry such that
$l$ meets the geodesic curves $g_{A_1A_2}^X$ in $Q$, $g_{A_0A_2}^X$ in $R$,
and $g_{A_0A_1}^X$ in $P$,
then $$s^X_g(A_0,P,A_1)s^X_g(A_1,Q,A_2)s^X_g(A_2,R,A_0) = -1.$$
\end{Theorem}
I would like to note here that the above theorems in fibre like cases are the same as the Euclidean case, are not mentioned here (see \cite{Sz22-2}).
\end{enumerate}
\begin{rmrk}
It is easy to see that the ``reversals'' of the above theorems are also true.
\end{rmrk}
\begin{rmrk}
The Ceva's and Menelaus' theorems for the {\it geodesic triangles} of the $\NIL$ space were examined in the article \cite{Sz23-1}.
\end{rmrk}
\subsection{Generalizations of Menelaus' and Ceva's theorems in $\NIL$ geometry}
We consider a {\it translation triangle}
$A_0A_1A_2$ in the projective model of the $\NIL$ space (see Section 2.1).
Without limiting generality, we can assume that $A_0=(1,0,0,0)$.
The translation curves that contain the sides $A_0A_1$ and $A_0A_2$ of the given triangle can be characterized directly by the corresponding parameters $\theta$ 
and $\phi$ (see (2.8) and (2.9)).

The translation curve including the translation-like side segment $A_1A_2$ is also determined by one of its endpoints and its parameters 
but in order to determine the corresponding parameters of this
translation curve we use the translation  $\bT^{\NIL}(A_1)$, as elements of the isometry group of the geometry $\NIL$, that
maps $A_1=(1,x_1,y_1,z_1)$ onto $A_0=(1,0,0,0)$. 

\begin{rmrk}
I note here, that we can use other orientation preserving isometry, but in all cases the $\theta$ parameters of the
``image translation curves" are equal (of course the $\phi$ parameters may be different) and in the further derivation, only the values of the parameter $\theta$ will be needed.
\end{rmrk}
From the (2.8), (2.9) equations of the $\NIL$ translation curve we directly obtain the following
\begin{lemma}
\begin{enumerate}
\item Let $P_1$ and $P_2$ be arbitrary points and $g_{P_1,P_2}^\NIL$ is the corresponding translation curve in the considered model of $\NIL$ geometry. 
The projection of the translation curve segment $g_{P_1,P_2}^\NIL$ onto the $[x,y]$ coordinate plane 
in the direction of the fibers will be a segment in the Euclidean sense if the points of the translation curve do not lie on a fiber
otherwise, its projection is a point (see Fig.~1).

\item If points $P_1$ and $P_2$ lie in the $[x,z]$ or $[y,z]$ coordinate plane, then the corresponding translation curve
$g_{P_1,P_2}^{\NIL}$ is a straight line segment $P_1P_2$ in the Euclidean sense (see Fig.~2).
\end{enumerate}~ ~ $\square$
\end{lemma}
\begin{rmrk}
In other words, the points of the translation curve $g_{P_1,P_2}^\NIL$ lie exactly in a plane in Euclidean sense orthogonal to the $[x,y]$ coordinate plane
if the points of the translation curve do not lie on a fiber.
\end{rmrk}
\begin{figure}[ht]
\centering
\includegraphics[width=10cm]{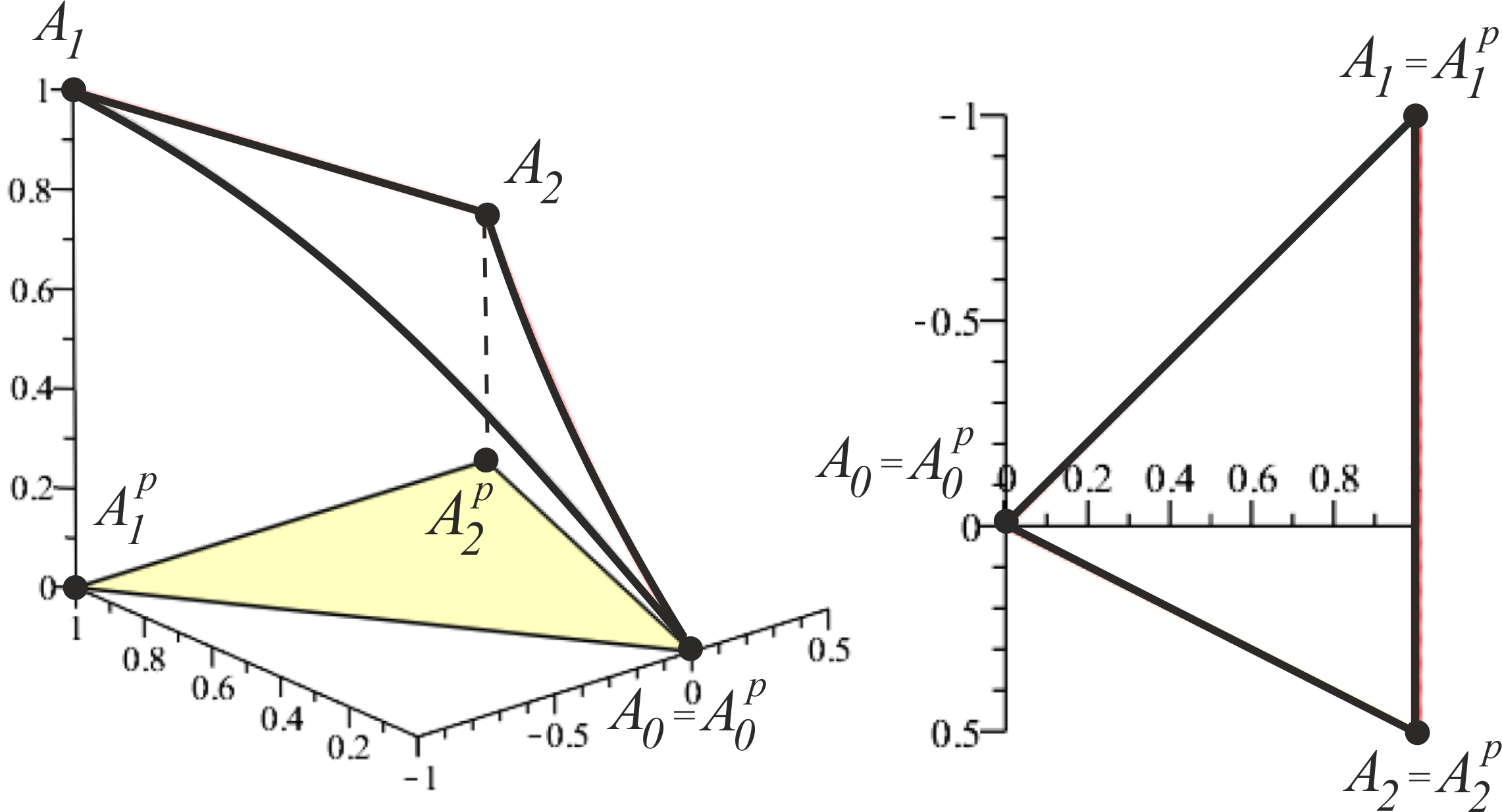}
\caption{Translation triangle with vertices $A_0=(1,0,0,0)$, $A_1=(1,-1,1,1)$, $A_2=(1,1/2,1,1/2)$ and its projected image by fibers into $[x,y]$ coordinate plane.}
\label{}
\end{figure}
\begin{figure}[ht]
\centering
\includegraphics[width=10cm]{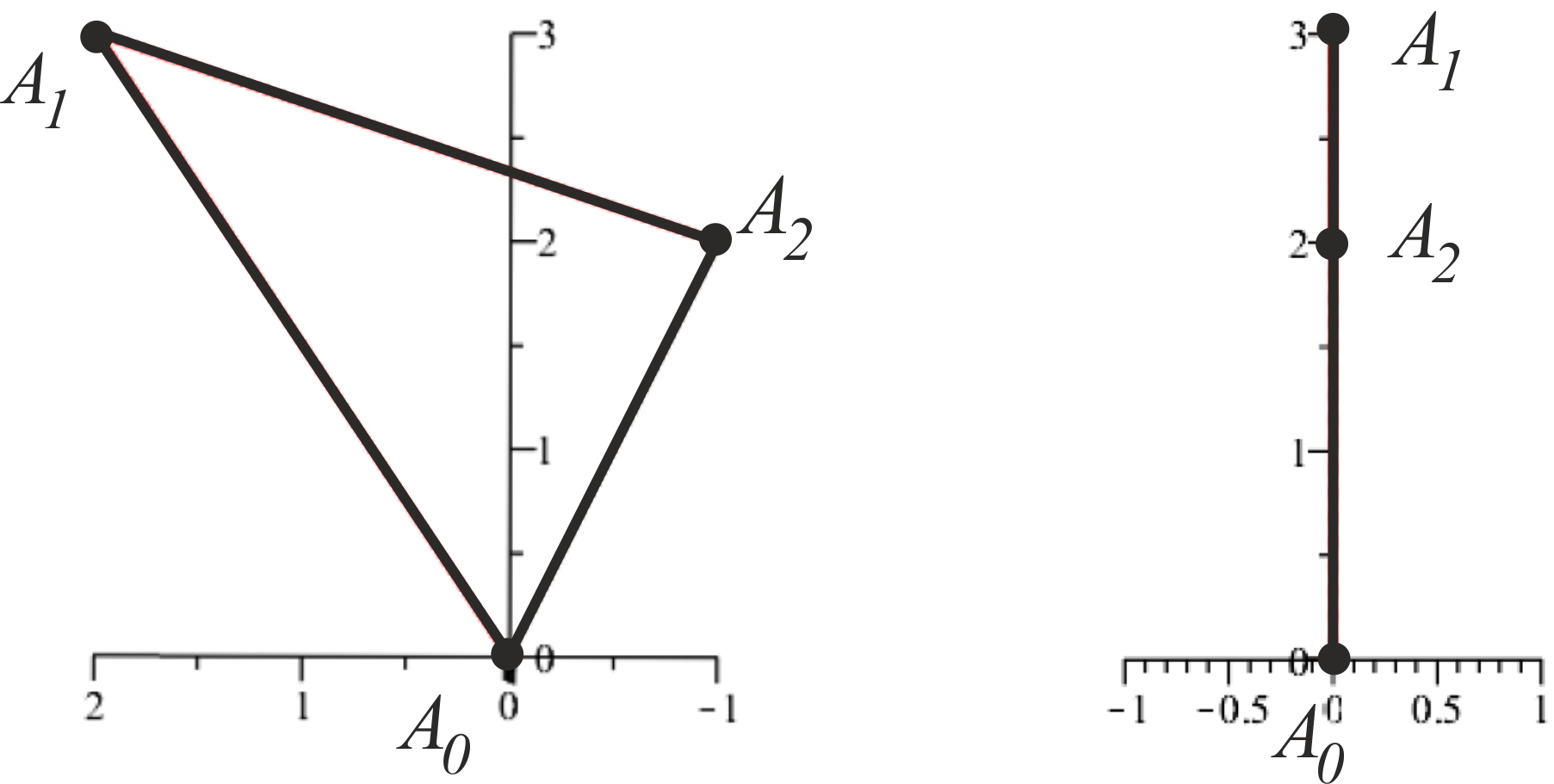}
\caption{Translation triangle with vertices $A_0=(1,0,0,0)$, $A_1=(1,2,0,3)$, $A_2=(1,-1,0,2)$ lying in the $[x,z]$ coordinate plane.}
\label{}
\end{figure}
We extend the definition of the simple ratio to the $\NIL$ space for translation curves.
\begin{definition}
If $A$, $B$ and $P$ are distinct points on a translation curve in the $\NIL$ space, then
their simple ratio is
$$s_g^\NIL(A,P,B) =  d^{\NIL,t}(A,P)/d^{\NIL,t}(P,B),$$ if $P$ is between $A$ and $B$, and
$$s_g^\NIL(A,P,B) = -d^{\NIL,t}(A,P)/d^{\NIL,t}(P,B),$$ (see Fig.~3).
\end{definition}
\begin{figure}[ht]
\centering
\includegraphics[width=7cm]{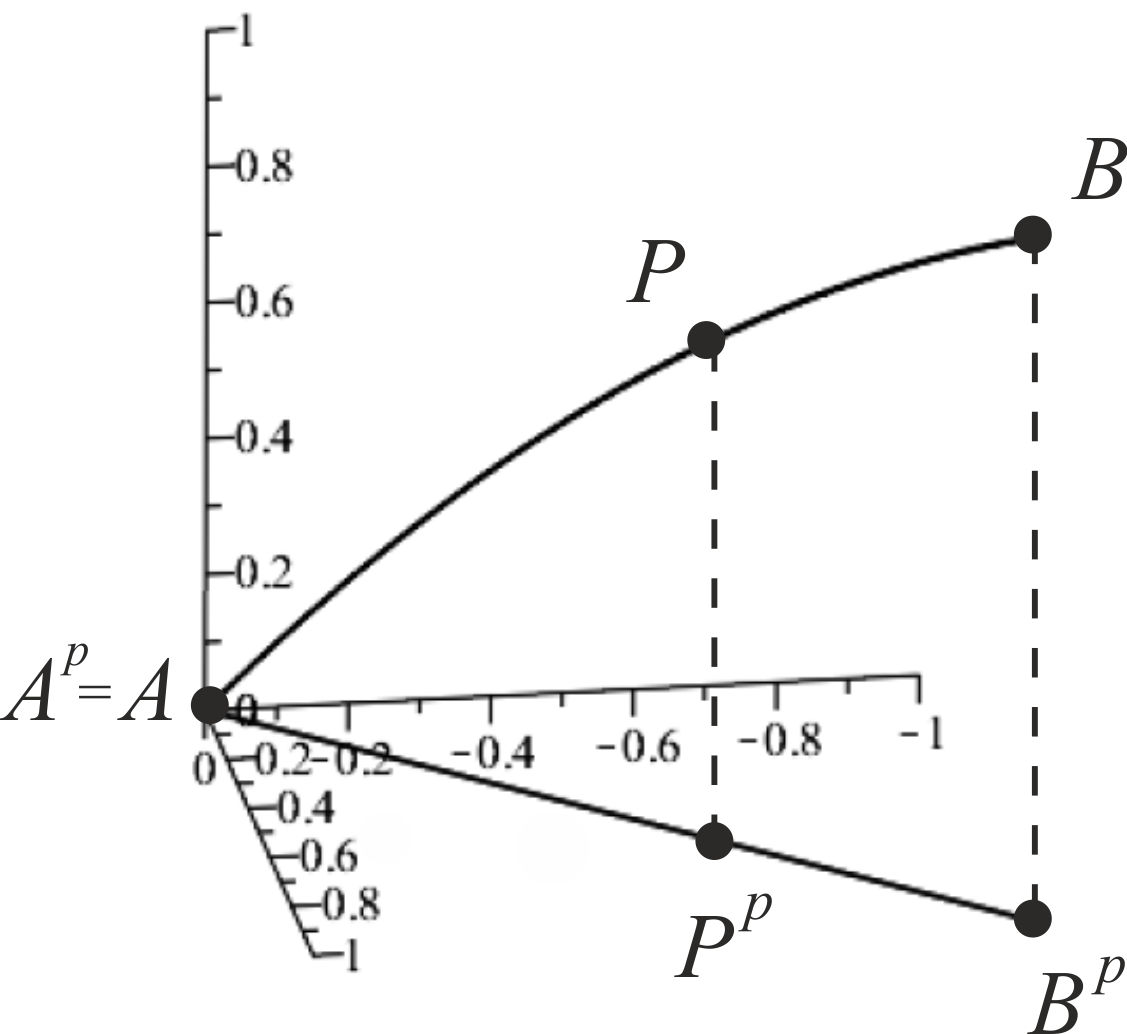}
\caption{$\NIL$ translation curve segment $g^\NIL_{A,B}$ with point $P \in g^\NIL_{A,B}$ and its projected image by fibers into $[x,y]$ coordinate plane.}
\label{}
\end{figure}
\begin{lemma}
\begin{enumerate}
\item Let $g_{A,B}$ be an arbitrary non-fibrum-like translation curve, where without loss of generality we can assume 
that $A=(1,0,0,0)$ coincides with the center of the model. 
Furthermore, let $P \in g^\NIL_{A,B}$ and let $A$, $B$ be the projected images 
by fibers into the $[x,y]$ coordinate plane, $A^p$, $B^p$, and $P^p$ then
$$
s_g^\NIL(A,P,B)=s_g^\NIL(A^p,P^p,B^p)=s_g^{\EUC}(A^p,P^p,B^p),
$$ 
{(see Definitions 3.13 and Lemma 3.11.)}
\item If points $A$, $B$ and $P$ lie on a fibrum, then
$$
s_g^\NIL(A,P,B)=s_g^{\EUC}(A,P,B)
$$
\end{enumerate}
\end{lemma}
{\bf Proof:}
\begin{enumerate}
\item From the Lemma (3.11) it is seen that the points $A^p$, $B^p$ and $P^p$ are located on a straight line in the $[x,y]$ plane. 
From the formulas (2.8), (2.9) it is directly seen that the translation distances on straight line $A^pB^p$ are equal to the Euclidean 
distances and that the ratio of the translation distances $d^t(A,P)$ and $d^t(P,B)$ is equal 
to the ratio of the translation distances $d^t(A^p,P^p)$ and $d^t(P^p,B^p)$. 
\item It follows directly from formulas (2.8) and (2.9), since the Euclidean distance on fibrum lines is equal to the translation distance.
\end{enumerate} \quad \quad $\square$

Let us introduce the following notations:

1. If the surface of a translation-like triangle is a plane in Euclidean sense
then it is called fibre type triangle.

2. In the other cases the triangle is in general type (see \cite{Cs-Sz25}).
\begin{definition}
Let $\mathcal{S}^{\NIL,t}_{A_0A_1A_2}$ be the surface of the translation triangle $A_0A_1A_2$ and $P_1$, $P_2 \in \mathcal{S}^{\NIL,t}_{A_0A_1A_2}$ two given point.

1. If the translation triangle is in fibre type, then the connecting curve
$P_1P_2 \subset \mathcal{S}^{\NIL,t}_{A_0A_1A_2}$ is a translation curve in $\NIL$ space.

2. In other cases the connecting curve $P_1P_2$ is the image of the translation curve $g^\NIL_{P_1P_2}$ into the surface 
$\mathcal{S}^{\NIL,t}_{A_0A_1A_2}$ by {\it fibrum projection}.
\end{definition}
\begin{Theorem}[Ceva's theorem for translation triangles in $\NIL$ space]
If $T$ is a point not on any side of a translation triangle $A_0A_1A_2$ in $\NIL$ such that
the curves $A_0T$ and $g^\NIL_{A_1,A_2}$ meet in $Q$, $A_1T$ and $g_{A_0A_2}^\NIL$ in $R$, and $A_2T$ and $g_{A_0A_1}^\NIL$ in $P$, $(A_0T, A_1T, A_2T \subset 
\mathcal{S}^\NIL_{A_0A_1A_2})$ 
then $$s_g^\NIL(A_0,P,A_1)s_g^\NIL(A_1,Q,A_2)s_g^\NIL(A_2,R,A_0) = 1.$$
\end{Theorem}
{\bf{Proof:}}

\begin{enumerate}
\item First we consider the translation triangle $A_0A_1A_2$ {\it triangle in general type} 
and let $A_0^pA_1^pA_2^p$ $(A_0^p=A_0)$ its fibrum-like projected image into the $[x,y]$ coordinate plane (see Fig.~1).

Moreover, let $T$ is a point not on any side of the translation triangle $A_0A_1A_2$ such that
the curves $A_0T$ and $g^\NIL_{A_1,A_2}$ meet in $Q$, $A_1T$ and $g_{A_0A_2}^\NIL$ in $R$, and $A_2T$ and $g_{A_0A_1}^\NIL$ in $P$, $(A_0T, A_1T, A_2T \subset 
\mathcal{S}^\NIL_{A_0A_1A_2})$ 

By the Lemma 3.11 and the Definition 3.15 their projected images $A_0^p Q^p$, $A_1^pR^p$ and $A_2^pP^p$ of the curves $A_0Q$, $A_1R$ and $A_2P$ are translation curves (straight segments) 
in the $[x,y]$ coordinate plane similarly to the side segments $A_0^pA_1^p$, $A_1^pA_2^p$, $A_2^pA_0^p$. Moreover,
from the Lemma 3.14 follows, that this projection does not change the simple ratio. 

Thus, the theorem follows directly from the Euclidean Ceva theorem. ~ ~ 
\item If the triangle is of fibrum type, then we get directly an Euclidean-like case (see Fig.~2), in which, by applying the Lemmas 3.11 and 3.14, obtain the theorem directly.
\quad \quad $\square$
\end{enumerate}
\begin{Theorem}[Menelaus's theorem for translation triangles in $\NIL$ space]
If $l$ is a line not through any vertex of a translation triangle
$A_0A_1A_2$ lying in a surface $\mathcal{S}^{\NIL,t}_{A_0A_1A_2}$
such that
$l$ meets the translation curves $g_{A_1A_2}^X$ in $Q$, $g_{A_0A_2}^X$ in $R$,
and $g_{A_0A_1}^\NIL$ in $P$,
then $$s^\NIL_g(A_0,P,A_1)s^\NIL_g(A_1,Q,A_22)s^\NIL_g(A_2,R,A_0) = -1.$$
\end{Theorem}
{\bf{Proof:}} 

Similarly to the above proof, this theorem is follows by the Definition 3.13 of the simple ratio , Lemmas 3.11 and 3.14 and the corresponding
Euclidean Menelaus' theorem. ~ ~ $\square$
\begin{rmrk}
It is easy to see that the ``reversals'' of the above theorems are also true.
\end{rmrk}
\subsection{Generalizations of Menelaus' and Ceva's theorems in $\SOL$ geometry}
We consider a {\it translation triangle}
$A_0A_1A_2$ in the projective model of the $\SOL$ space (see Section 2.2).
Without limiting generality, we can assume that $A_0=(1,0,0,0)$.
The translation curves that contain the sides $A_0A_1$ and $A_0A_2$ of the given triangle can be characterized directly by the corresponding parameters $\theta$ 
and $\phi$ (see (2.20)).

The translation curve including the translation-like side segment $A_1A_2$ is also determined by one of its endpoints and its parameters 
but in order to determine the corresponding parameters of this
translation line we use the translation  $\bT^{\SOL}(A_1)$, as elements of the isometry group of the geometry $\SOL$, that
maps $A_1=(1,x_1,y_1,z_1)$ onto $A_0=(1,0,0,0)$. 

From the equation system (see (2.20)) of the $\SOL$ translation curve we directly obtain the following
\begin{lemma}
\begin{enumerate}
\item Let $A_0=(1,0,0,0)$ and $A_1=(1,x_1,y_1,z_1)$ be points ($A_1 \notin [y,z]$ or $[x,y]$ coordinate planes) 
and $g_{A_0,A_1}^\SOL$ is the corresponding translation curve in the considered model of $\SOL$ geometry (see Section 2.2, (2.20)). 
Let the orthogonal projected image in Euclidean sense of the translation curve segment $g_{A_0,A_1}^\SOL$ onto the $[x,z]$ coordinate plane be $g_{A_0,A_1}^{p,\SOL}$.
Then the equation of $g_{A_0,A_1}^{p,\SOL}$ in the $[x,z]$ coordinate plane is:
\begin{equation}
z=-\log\Big(\frac{-x}{\cos \theta \cos \phi}+1\Big), 
x \in (-\cot \theta \cos \phi,\cot \theta \cos \phi), \tag{3.1}
\end{equation}
where $\theta$ and $\phi$ are the parameters of $g_{A_0,A_1}^\SOL$ (see 2.20). 
\item If the points $A_0$ and $A_1$ lie in the $[y,z]$ coordinate plane, then the corresponding translation curve
$g_{A_0,A_1}^\SOL$ lies in the coordinate plane $[y,z]$ and its equation is:
\begin{equation}
z=\log\Big(\frac{y}{\cos \theta \sin \phi}+1\Big), 
y \in (-\cot \theta \sin \phi,\cot \theta \sin \phi),  \tag{3.2}
\end{equation}
where $\theta$ and $\phi=\pi/2$ are the parameters of $g_{A_0,A_1}^\SOL$ (see 2.20). 
\item If the points $A_0$ and $A_1$ lie in the $[x,y]$ coordinate plane, then the corresponding translation curve
$g_{A_0,A_1}^\SOL$ lies in the coordinate plane $[x,y]$ and its equation is:
\begin{equation}
y=\tan\phi \cdot x, ~ x \in \bR, \tag{3.3}
\end{equation}
where $\theta=0$ and $\phi$ are the parameters of $g_{A_0,A_1}^\SOL$ (see 2.20). 
\end{enumerate}~ ~ $\square$
\end{lemma}
\begin{figure}[ht]
\centering
\includegraphics[width=12cm]{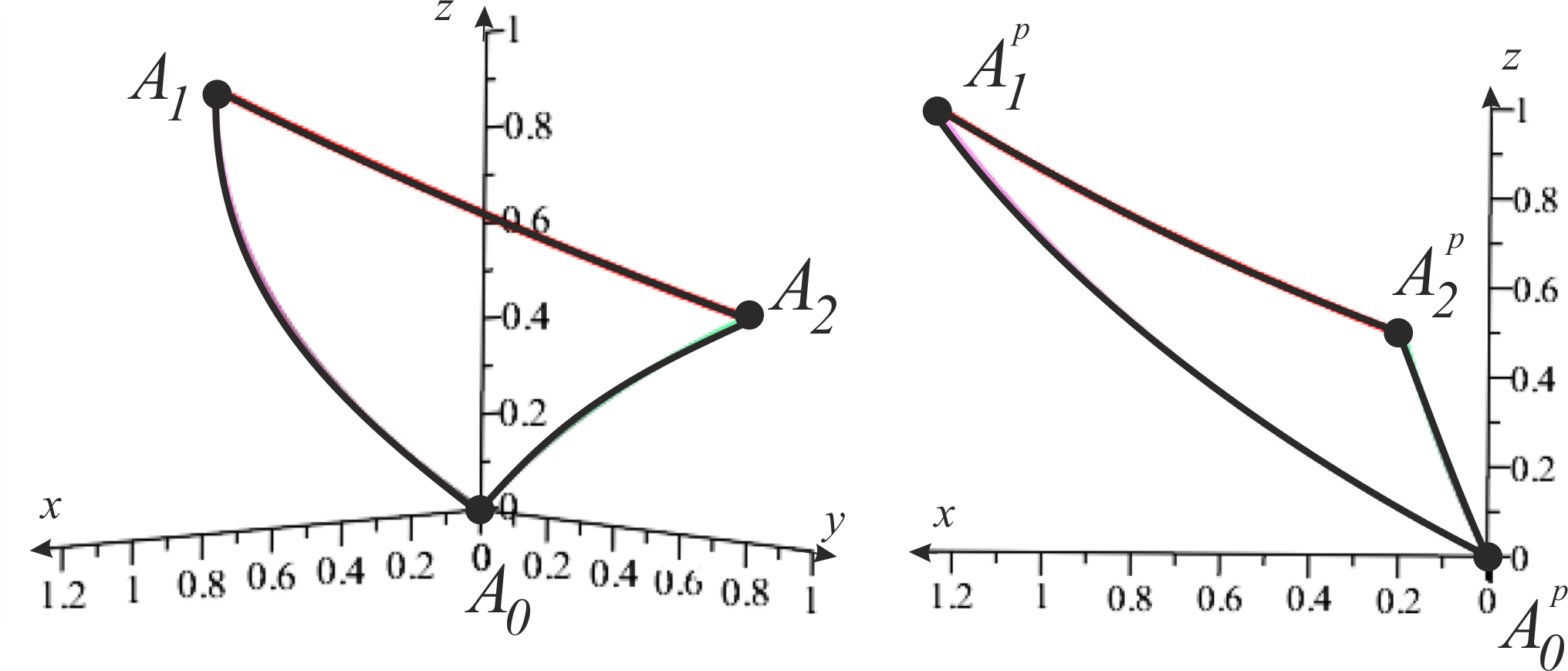}
\caption{Translation triangle with vertices $A_0=(1,0,0,0)$, $A_1=(1,5/4,1/2,1)$, $A_2=(1,1/5,1,1/2)$ and its image by orthogonal projection into $[x,z]$ coordinate plane.}
\label{}
\end{figure}
\begin{figure}[ht]
\centering
\includegraphics[width=12cm]{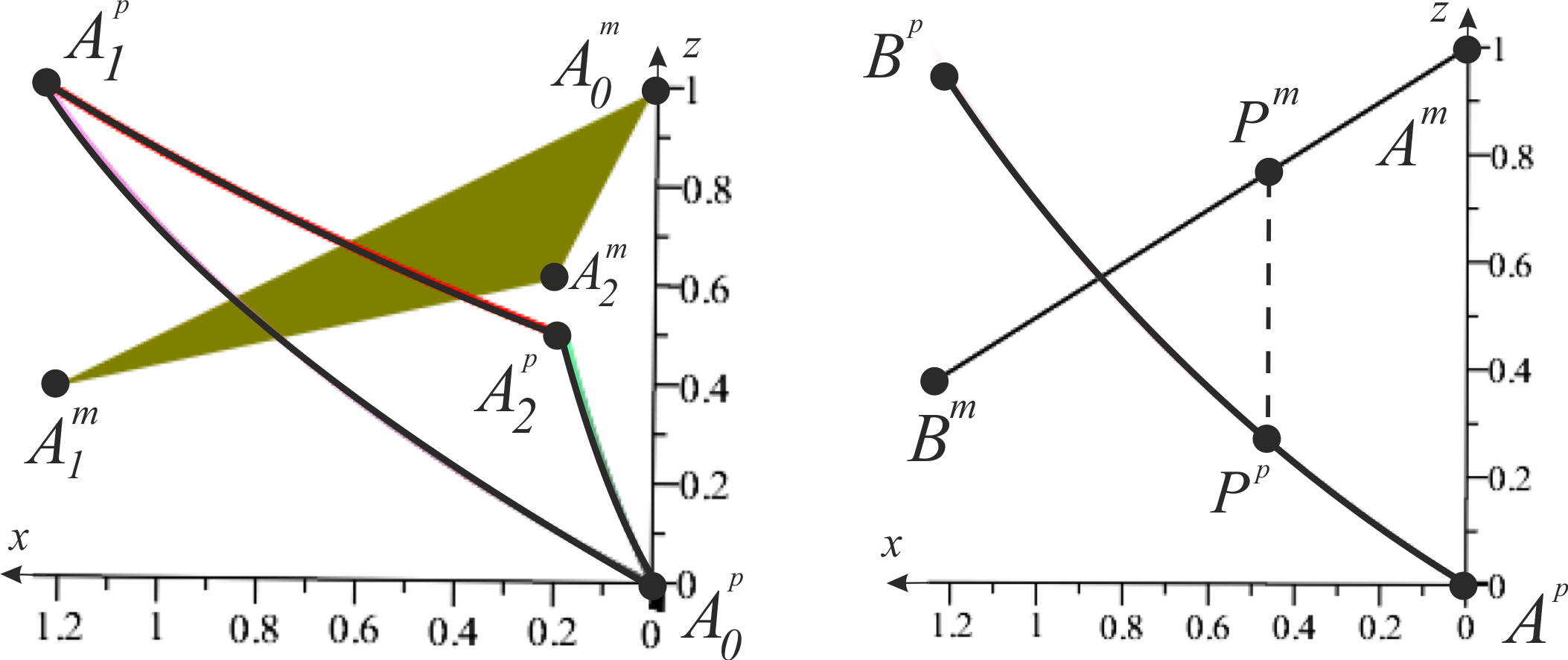}
\caption{The projected translation triangle with vertices $A_0^p=(1,0,0,0)$, $A_1^p=(1,5/4,0,1)$, $A_2^p=(1,1/5,0,1/2)$ lying in the $[x,z]$ coordinate plane and 
its transformed image with vertices
$A_0^m=(1,0,0,1)$, $A_1^m=(1,5/4,0,e^1)$, $A_2^m=(1,1/5,0,e^{1/2})$ (left). The projected $\SOL$ translation curve 
segment $g^{p,\SOL}_{A,B}$ with point $P^p \in g^{p,\SOL}_{A,B}$ and its transformed imege with points $A^m, B^m, P^m$ in $[x,z]$ coordinate plane (right).}
\label{}
\end{figure}

The metric on the plane $[x,z]$ induced by $(ds)^2 =e^{2z}(dx)^2 +(dz)^2$ agrees with the $\SOL$ metric.
On mapping this plane to the upper half plane $\{(x_1,x_2) \in \mathbf{R}:~x_2 > 0$ by mapping 
$ma:(x,z) \mapsto (x_1,x_2) =(x,e^{-z})$. We obtain, that the metric becomes 
the standard hyperbolic metric. 
\begin{equation}
(ds)^2 = \frac{(dx_1)^2+(dx_2)^2}{x_2^2}. \notag
\end{equation}
Therefore, the $[x,z]$ plane (and so the plane $[y,z]$) is a convexly embedded copy 
of $\mathbf{H}^2$ in $\SOL$. Moreover, a similar statement can be formulated for the planes $x=c$, $y=c$ $(c\in \mathbf{R})$.
Moreover, by Lemma 3.19 we directly obtain the following
\begin{lemma}
\begin{enumerate}
\item Let $A_0=(1,0,0,0)$ and $A_1=(1,x_1,y_1,z_1)$ be points ($A_1 \notin y$ axis) 
and $g_{A_0,A_1}^\SOL$ is the corresponding translation curve in the considered model of $\SOL$ geometry (see Section 2.2, (2.20)). 
Let the orthogonal projected image in Euclidean sense of the translation curve segment $g_{A_0,A_1}^\SOL$ onto the $[x,z]$ coordinate plane be $g_{A_0,A_1}^{p,\SOL}$.
If we apply the mapping $m$ to the points of $g_{A_0,A_1}^{p,\SOL}$ then its image is a segment $A_0^m A_1^m$ in the Euclidean sense in the coordinate plane $[x,z]$ (see Fig.~5)
with equation
$$
z=\frac{-x}{\cos \theta \cos \phi}+1. 
$$
where $\theta$ and $\phi$ are the parameters of $g_{A_0,A_1}^\SOL$ (see 2.20). 
\item If the points $A_1 \in y$ ($A_0A_1$ is an straight segment in Euclidean sense) then its orthogonal projected image onto the $[x,z]$ coordinate plane is a point and therefore its 
the image under transformation $m$ will also be a point with coordinates $(1,0,0,1)$.
\end{enumerate}~ ~ $\square$
\end{lemma}
\begin{definition}
Let $\mathcal{S}^{\SOL,t}_{A_0A_1A_2}$ be the surface of the translation triangle $A_0A_1A_2$ (can be assumed that $A_0=(1,0,0,0))$ and $P_1$, 
$P_2 \in \mathcal{S}^{\SOL,t}_{A_0A_1A_2}$ two given point.
\begin{enumerate}
\item If the vertices $A_1, A_2$ of the translation triangle $A_0A_1A_2$ lie on the $[y,z]$ or $[x,z]$ coordinate plane, their connecting curve
$P_1P_2 \subset \mathcal{S}^{\SOL,t}_{A_0A_1A_2}$ is a translation curve in $\SOL$ space that lies on the corresponding coordinate plane.

\item In other cases the connecting curve $P_1P_2$ is derived by the following steps:
\begin{enumerate}
\item Let the perpendicular projections of the points $P_1$ and $P_2$ onto the $[x,z]$ plane be $P_1^p$ and $P_2^p$.
If $P_1^p \ne P_2^p$ then apply the mapping $m$ to these points: 
$m: P_1^p \rightarrow P_1^m$, $m: P_2^p \rightarrow P_2^m$. 
\item
Denote the Euclidean line of $P_1^m P_2^m$ by $g_{P_1,P_2}^{m,\SOL}$ and its image at $m^{-1}$ by $g_{P_1,P_2}^{p,\SOL}$.
\item
Let $\cC \cS_{P_1,P_2}$ be a {\it cylindrical surface} consisting of all the points on all the lines which are parallel to $y$ axis
and pass through the points of the plane curve $g_{P_1,P_2}^{p,\SOL}$. 
\item The connecting curve $P_1P_2 \subset \mathcal{S}^{\SOL,t}_{A_0A_1A_2}$ is 
$P_1P_2:=\mathcal{S}^{\SOL,t}_{A_0A_1A_2} \cap \cC \cS_{P_1,P_2}.$
\end{enumerate}
\item If $P_1^p = P_2^p$ then let the curve connecting $P_1^p$ and $P_2^p$ (similarly $P_1^m$ and $P_2^m$) be parallel to the $z$ axis.
Let $\cC \cS_{P_1,P_2}$ be a Euclidean plane parallel to $[y,z]$ coordinate plane which contains points $P_1^p$ and $P_2^p$. 
			       
The connecting curve $P_1P_2 \subset \mathcal{S}^{\SOL,t}_{A_0A_1A_2}$ is 
$P_1P_2:=\mathcal{S}^{\SOL,t}_{A_0A_1A_2} \cap \cC \cS_{P_1,P_2}.$
\end{enumerate}
\end{definition}
We extend the definition of the simple ratio to the $\SOL$ space for translation curves. 
\begin{definition}
\begin{enumerate}
\item If $A$, $B$ and $P$ are distinct points on a translation curve in the $\SOL$ space where this translation curve does not lie in the coordinate plane $[x,y]$ (see Fig.~5), then
their simple ratio is
$$s_g^\SOL(A,P,B) =  \frac{1-e^{-d^{\SOL,t}(A,P)\cdot\sin{\theta}}}{e^{-d^{\SOL,t}(A,P)\cdot\sin{\theta}}-e^{-d^{\SOL,t}(A,B)\cdot\sin{\theta}}},$$ if $P$ is between $A$ and $B$, and
$$s_g^\SOL(A,P,B) = -\frac{1-e^{-d^{\SOL,t}(A,P)\cdot\sin{\theta}}}{e^{-d^{\SOL,t}(A,P)\cdot\sin{\theta}}-e^{-d^{\SOL,t}(A,B)\cdot\sin{\theta}}},$$ 
otherwise where
$\theta$ is
the parameter of the translation curve containing points $A,B,P$.
\item If $A$, $B$ and $P$ are distinct points on a translation curve in the $\SOL$ space where this translation curve lies in the coordinate plane $[x,y]$ or in the $z$-axis 
(see (2.20) then
$$s_g^\SOL(A,P,B) = d^{\SOL,t}(A,P)/d^{\SOL,t}(P,B),$$ if $P$ is between $A$ and $B$, and
$$s_g^\SOL(A,P,B) = -d^{\SOL,t}(A,P)/d^{\SOL,t}(P,B)$$.
\item Let $A$, $B$ and $P$ be distinct points on a translation curve in the $\SOL$ space where this translation curve lies in the coordinate plane $[y,z]$. 
For this translation curve we can apply an isometry of $\SOL$ space that assigns the $[x,z]$ plane to the $[y,z]$ plane (see (2.13), (2.14)) 
therefore the simple ratio is derived by 
the first point of this definition.
\end{enumerate}
\end{definition}
\begin{lemma}
\begin{enumerate}
\item Let $A=(1,0,0,0)$ and $B=(1,x_1,y_1,z_1)$ be points ($B \notin [y,z]$ coordinate plane) 
and $g_{A,B}^{\SOL}$ is the corresponding translation curve in the considered model of $\SOL$ geometry (see Section 2.2, (2.20)). 
Moreover, let the orthogonal projected image in Euclidean sense of the translation curve segment $g_{A,B}^\SOL$ onto the $[x,z]$ coordinate plane be $g_{A,B}^{p,\SOL}$.
If we apply the mapping $m$ to the points of $g_{A,B}^{p,\SOL}$ then we obtain the segment $A^m B^m$ in the Euclidean sense in the coordinate plane $[x,z]$ 
(see Fig.~5) then we obtain the followig:
$$
s_g^\SOL(A,P,B)=s_g^{\EUC}(A^m,P^m,B^m),
$$ 
{(see Definitions 3.21 and 3.22)}.
\item If points $A_0$, $A_1$ and $P$ lie on $[x,y]$ coordinate plane or in $z$ axis, then
$$
s_g^\SOL(A,P,B)=s_g^{\EUC}(A,P,B).
$$
\item
Let $A$, $B$ and $P$ be distinct points on a translation curve in the $\SOL$ space where this translation curve lies in the coordinate plane $[y,z]$. 
For this translation curve we can apply an isometry of $\SOL$ space that assigns the $[x,z]$ plane to the $[y,z]$ plane (see (2.13), (2.14)) 
therefore similarly to the first point of the lemma, the equation given there is satisfied.
\end{enumerate}
\end{lemma}
{\bf Proof:}
\begin{enumerate}
\item Let the orthogonal projected image in Euclidean sense of the points 
$A^p,B^p,P^p$ onto the $x$ axis be $A^{p,x},B^{p,x},P^{p,x}$. 
It is clear, that $s_g^{\EUC}(A^m,P^m,B^m)=s_g^{\EUC}(A^{p,x},B^{p,x},P^{p,x})$. 
From the formula (2.20) and by Lemma 3.19 it is directly seen that the 
translation distances on $x$ axis are equal to the Euclidean 
distances and thus the simple ratio $s_g^{\EUC}(A^{p,x},B^{p,x},P^{p,x})
=s_g^{\SOL}(A^{p,x},B^{p,x},P^{p,x})$. 
From this, the statement of the theorem follows by direct 
application of formula (3.1) (and similarly the 
connection between the translation distances $d^\SOL(A,P)$,  
$d^\SOL(P,B)$ and the translation (or Euclidean) distances 
$d^{\SOL,t}(A^{p,x},P^{p,x})$, $d^{\SOL,t}(P^{p,x},B^{p,x})$).
\item It follows directly from formulas (2.20), since the Euclidean distance on $[x,y]$ plane and in $z$ axis is equal to the translation distance.
\end{enumerate} \quad \quad $\square$
\begin{Theorem}[Ceva's theorem for translation triangles in $\SOL$ space]
If $T$ is a point not on any side of a translation triangle $A_0A_1A_2$ in $\SOL$ 
such that the curves $A_0T$ and $g^\SOL{A_1,A_2}$ meet in $Q$, $A_1T$ and 
$g_{A_0A_2}^\SOL$ in $R$, and $A_2T$ and $g_{A_0A_1}^\SOL$ in $P$, 
$(A_0T, A_1T, A_2T \subset 
\mathcal{S}^\SOL_{A_0A_1A_2})$ 
then $$s_g^\SOL(A_0,P,A_1)s_g^\SOL(A_1,Q,A_2)s_g^\SOL(A_2,R,A_0) = 1.$$
\end{Theorem}
{\bf{Proof:}}
\begin{enumerate}
\item We consider a {\it translation triangle}
$A_0A_1A_2$ in the projective model of the $\SOL$ space (see Section 2.2).
Without limiting generality, we can assume that $A_0=(1,0,0,0)$ and at most one of the points $A_1,A_2$ is on the $[y,z]$ coordinate plane. 

Moreover, let the orthogonal projected image in Euclidean sense of the translation curve segment $g_{A_i,A_j}^\SOL$ onto the $[x,z]$ coordinate plane be $g_{A_i,A_j}^{p,\SOL}$
$(i \ne j, i<j, i,j \in\{0,1,2\}$ and the images of points $T,Q,R,P$ are $T^p,Q^p,R^p,P^p$. 
If we apply the mapping $m$ to the above triangle and the corresponding points we obtain the Euclidean-like $A_1^m A_2^m A_3^m$ triangle and the points $T^m,Q^m,R^m,P^m$ 
lying in the $[x,z]$ coordinate plane. The Euclidean Ceva theorem is true for this configuration exactly if  
$${s_g}^{\EUC} (A_0^m,P^m,A_1^m) {s_g}^{\EUC} (A_1^m,Q^m,A_2^m) {s_g}^{\EUC} (A_2^m,R^m,A_0^m) = 1.$$ 

After this we applied to this arrangement, then the $m^{-1}$ transformation and 
we construct the $A_0Q, A_1R, A_2P \subset \mathcal{S}^{\SOL,t}_{A_0A_1A_2}$ surface curves as described in Definition 3.21, 
then we get that these curves on $\mathcal{S}^{\SOL,t}_{A_0A_1A_2}$ also pass through the point $T$.
This, due to the Lemma 3.23, occurs exactly when the definition of the simple ratio is based on Definition 3.22. 
\item If the triangle lies on the coordinate plane $[x,y]$, then we get directly an Euclidean-like case, in which, by applying the Lemmas 3.19 and 3.23, obtain the theorem directly.
\quad \quad $\square$
\end{enumerate}
\begin{Theorem}[Menelaus's theorem for translation triangles in $\SOL$ space]
If $l$ is a line not through any vertex of a geodesic triangle
$A_0A_1A_2$ lying in a surface $\mathcal{S}^{\SOL,t}_{A_0A_1A_2}$
such that
$l$ meets the geodesic curves $g_{A_1A_2}^\SOL$ in $Q$, $g_{A_0A_2}^SOL$ in $R$,
and $g_{A_0A_1}^\SOL$ in $P$,
then $$s^\SOL_g(A_0,P,A_1)s^\SOL_g(A_1,Q,A_22)s^\SOL_g(A_2,R,A_0) = -1.$$
\end{Theorem}
{\bf{Proof:}} 

Similarly to the above proof, this theorem is follows by the Definition 3.21 of the simple ratio , Lemmas 3.19 and 3.23 and the corresponding
Euclidean Menelaus' theorem. ~ ~ $\square$
\begin{rmrk}
It is easy to see that the ``reversals'' of the above theorems are also true.
\end{rmrk}
\subsection{Generalizations of Menelaus' and Ceva's theorems in $\SLR$ geometry}
We consider a {\it translation triangle}
$A_0A_1A_2$ in the projective model of the $\SLR$ space (see Section 2.3).
Without limiting generality, we can assume that $A_0=(1,0,0,0)$.
The translation curves that contain the sides $A_0A_1$ and $A_0A_2$ of the given triangle can be characterized directly by the corresponding parameters $\alpha$ 
and $\lambda$ (see Table 1).

The translation curve including the translation-like side segment $A_1A_2$ is also determined by one of its endpoints and its parameters 
but in order to determine the corresponding parameters of this translation curve we use the translation  $\bT^{\SLR}(A_1)$, as elements of the isometry group of the geometry $\SOL$, that
maps $A_1=(1,x_1,y_1,z_1)$ onto $A_0=(1,0,0,0)$. From the equation of the translation curves, 
it can be seen that they are straight lines in the Euclidean sense, so the translation triangle will also be a triangle in the 
Euclidean sense in our model. Moreover, the surface of the translation triangle will also be a plane in the Euclidean sense (see \cite{Cs-Sz24}).

Thus, we can generalize the Ceva's and Menelaus' theorems by directly extending the corresponding Euclidean theorems. To do this, we 
need to determine the relationship between the concepts of Euclidean and $\SLR$ simple relations. 

It can be assumed by the homogeneity of $\SLR$ that the starting point of a 
a given translation curve segment is $A=(1,0,0,0)$ and 
the other endpoint will be given by its homogeneous coordinates $B=(1,a,b,c)$. 
We consider the translation curve segment $g_{A,B}^{\SLR}$ and determine its
parameters $\alpha,\lambda,s$ expressed by the real coordinates $a$, $b$, $c$ of $B$.We can assume by the structure of $\SLR$ that $a \in \bR^+$.
We obtain directly by equation systems in Table 1. the following:
\begin{lemma}[\cite{Sz19}]
\begin{enumerate}
\item Let $(1,a,b,c)$ $(0 \ne b,c \in \bR, a \in \bR, a^2-b^2-c^2<0\})$ $(0 \le \alpha < \frac{\pi}{4}, ~ -\pi < \lambda \le \pi)$ be the homogeneous coordinates of the point $P \in \SLR$. The parameters of the
corresponding translation curve $t_{E_0P}$ are the following
\begin{equation}
\begin{gathered}
\alpha=\mathrm{arctg}\Big(\frac{a}{\sqrt{b^2+c^2}} \Big), \ \
\lambda=\mathrm{arctg}\Big(\frac{c}{b}\Big), \\
s=\frac{\mathrm{arctanh}\Big(\frac{a \cdot \sqrt{\cos{2\alpha}}}{\sin{\alpha}}\Big)}{\sqrt{\cos{2\alpha}}}=\\
=\frac{\mathrm{arctanh}\Big(\sqrt{a^2+b^2+c^2}\cdot \sqrt{\cos\Big(2 \cdot \mathrm{arctan} \Big(\frac{a}{\sqrt{b^2+c^2}}\Big)\Big)}}{\sqrt{\cos\Big(2\mathrm{arctan} \Big(\frac{a}{\sqrt{b^2+c^2}}\Big) \Big)}}
\end{gathered} \tag{3.4}
\end{equation}
\item Let $(1,a,b,c)$ $(0 \ne b,c \in \bR, a \in \bR^+, a^2-b^2-c^2=0\})$ $-\pi < \lambda \le \pi)$ be the homogeneous coordinates of the point $P \in \SLR$. The parameters of the
corresponding translation curve $t_{E_0P}$ are the following
\begin{equation}
\begin{gathered}
\alpha=\frac{\pi}{4}, \ \ \lambda=\mathrm{arctg}\Big(\frac{c}{b}\Big), \ \ \ s=\frac{a}{\sqrt{2}} ,\\
\end{gathered} \tag{3.5}
\end{equation}
\item Let $(1,a,b,c)$ $(0 \ne b,c \in \bR, a \in \bR^+, a^2-b^2-c^2>0\})$ $(\frac{\pi}{4} < \alpha < \frac{\pi}{2}, ~ -\pi < \lambda \le \pi)$ be the homogeneous coordinates of the point $P \in \SLR$. The parameters of the
corresponding translation curve $t_{E_0P}$ are the following
\begin{equation}
\begin{gathered}
\alpha=\mathrm{arctan}\Big(\frac{a}{\sqrt{b^2+c^2}} \Big), \ \
\lambda=\mathrm{arctan}\Big(\frac{c}{b}\Big), \\
s=\frac{\mathrm{arctan}\Big(\frac{a \cdot \sqrt{-\cos{2\alpha}}}{\sin{\alpha}}\Big)}{\sqrt{-\cos{2\alpha}}}=\\
=\frac{\mathrm{arctan}\Big(\sqrt{a^2+b^2+c^2}\cdot \sqrt{-\cos\Big(2 \cdot \mathrm{arctan} \Big(\frac{a}{\sqrt{b^2+c^2}}\Big)\Big)}}{\sqrt{-\cos\Big(2\mathrm{arctan} \Big(\frac{a}{\sqrt{b^2+c^2}}\Big) \Big)}}
\end{gathered} \tag{3.6}
\end{equation}
\item Let $(1,a,b,c)$ $(0 = b,c \in \bR, a \in \bR^+\})$ be the homogeneous coordinates of the point $P \in \SLR$. The parameters of the
corresponding translation curve $t_{E_0P}$ are the following
\begin{equation}
\begin{gathered}
\alpha=\frac{\pi}{2}, \ \ s=\mathrm{arctg}(x).
\end{gathered} \ \ \ \ \square \tag{3.7}
\end{equation}
\end{enumerate}
\label{invlem}
\end{lemma}
We extend the definition of the simple ratio to the $\SLR$ space for translation curves.
\begin{definition}
\begin{enumerate}
\item If the parameter the translation curve containing points $A,B,P$ is $\alpha \in [0,\pi/4)$ then 
their simple ratio is
$$s_g^{\SLR}(A,P,B) =  \frac{\mathrm{tanh}(d^{\SLR,t}(A,P)\cdot \sqrt{\cos{2\alpha}})}{\mathrm{tanh}(d^{\SLR,t}(P,B)\cdot 
\sqrt{\cos{2\alpha}})},$$ if $P$ is between $A$ and $B$, and
$$s_g^{\SLR}(A,P,B) =  -\frac{\mathrm{tanh}(d^{\SLR,t}(A,P)\cdot \sqrt{\cos{2\alpha}})}{\mathrm{tanh}(d^{\SLR,t}(P,B)\cdot 
\sqrt{\cos{2\alpha}})},$$ 
otherwise. 
\item If the parameter the translation curve containing points $A,B,P$ is $\alpha =\pi/4$ then 
their simple ratio is
$$s_g^{\SLR}(A,P,B) =  d^{\SLR,t}(A,P)/d^{\SLR,t}(P,B),$$ if $P$ is between $A$ and $B$, and
$$s_g^{\SLR}(A,P,B) =  -d^{\SLR,t}(A,P)/d^{\SLR,t}(P,B),$$ 
otherwise. 
\item If the parameter the translation curve containing points $A,B,P$ is $\alpha \in (\pi/4,\pi/2]$ then 
their simple ratio is
$$s_g^{\SLR}(A,P,B) =  \frac{\mathrm{tan}(d^{\SLR,t}(A,P)\cdot \sqrt{-\cos{2\alpha}})}{\mathrm{tan}(d^{\SLR,t}(P,B)\cdot 
\sqrt{-\cos{2\alpha}})},$$ if $P$ is between $A$ and $B$, and
$$s_g^{\SLR}(A,P,B) =  -\frac{\mathrm{tan}(d^{\SLR,t}(A,P)\cdot \sqrt{-\cos{2\alpha}})}{\mathrm{tan}(d^{\SLR,t}(P,B)\cdot 
\sqrt{-\cos{2\alpha}})},$$ 
otherwise. 
\end{enumerate}
\end{definition}
\begin{figure}[ht]
\centering
\includegraphics[width=10cm]{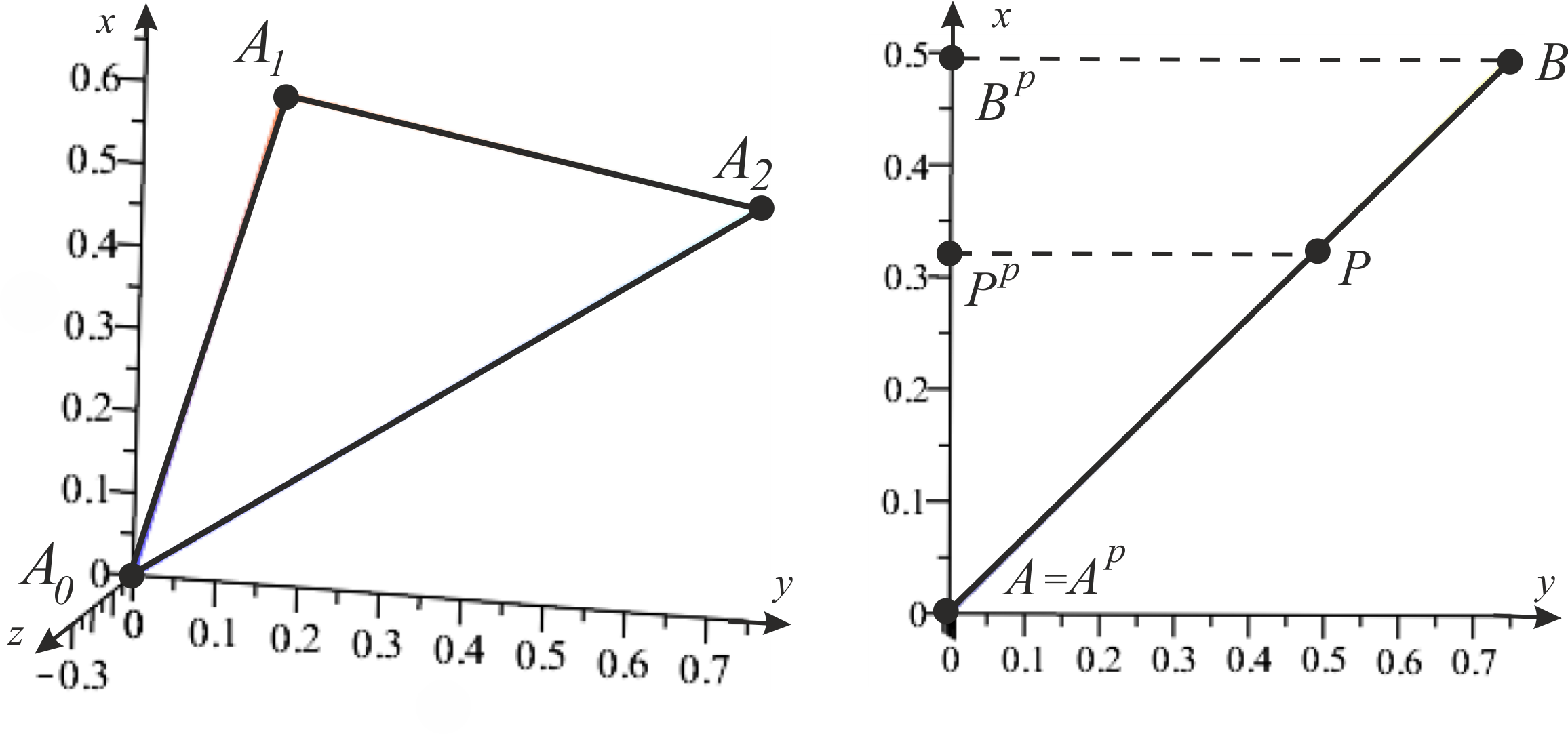}
\caption{Translation triangle with vertices $A_0=(1,0,0,0)$, $A_1=(1,1/2,3/4,0)$, $A_2=(1,2/3,1/4,-1/3)$ (left) The translation segment $g_{A,P,B}^{\SLR.t}$ and its image 
by orthogonal projection into $x$ axis.}
\label{}
\end{figure}
The following lemma is a direct consequence of Definition 3.28. and Lemma 3.27.:
\begin{lemma}
\begin{enumerate}
\item Let $g_{A,B}^{\SLR}$ be a translation curve which does not lie in the $[y,z]$ coordinate plane, where without loss of generality we can assume 
that $A=(1,0,0,0)$ coincides with the center of the model. 
Furthermore, let $P \in g^{\SLR}_{A,B}$ and let the orthogonally projected images of points $A,B,P$ 
into the $x$-axis be, $A^p$, $B^p$, and $P^p$ then
$$
s_g^{\SLR}(A,P,B)=s_g^{\EUC}(A^p,P^p,B^p),
$$ 
\item If points $A$, $B$ and $P$ lie on $[y,z]$ coordinate plane, then
$$
s_g^{\SLR}(A,P,B)=s_g^{\EUC}(A,P,B).
$$
\end{enumerate}
\end{lemma} \quad \quad $\square$

From the previous lemma, the Ceva's and Menelaus' theorems for translational triangles in $\SLR$ geometry can be directly derived.
\begin{Theorem}[Ceva's theorem for translation triangles in $\SLR$ space]
If $T$ is a point not on any side of a translation triangle $A_0A_1A_2$ in ${\SLR}$ 
such that the curves $A_0T$ and $g^\SOL{A_1,A_2}$ meet in $Q$, $A_1T$ and 
$g_{A_0A_2}^{\SLR}$ in $R$, and $A_2T$ and $g_{A_0A_1}^{\SLR}$ in $P$, 
$(A_0T, A_1T, A_2T \subset 
\mathcal{S}^{\SLR}_{A_0A_1A_2})$ 
then $$s_g^{\SLR}(A_0,P,A_1)s_g^{\SLR}(A_1,Q,A_2)s_g^{\SLR}(A_2,R,A_0) = 1.$$
\end{Theorem} \quad \quad $\square$
\begin{Theorem}[Menelaus's theorem for translation triangles in $\SOL$ space]
If $l$ is a line not through any vertex of a geodesic triangle
$A_0A_1A_2$ lying in a surface $\mathcal{S}^{{\SLR},t}_{A_0A_1A_2}$
such that
$l$ meets the geodesic curves $g_{A_1A_2}^{\SLR}$ in $Q$, $g_{A_0A_2}^{\SLR}$ in $R$,
and $g_{A_0A_1}^{\SLR}$ in $P$,
then $$s^{\SLR}_g(A_0,P,A_1)s^{\SLR}_g(A_1,Q,A_22)s^{\SLR}_g(A_2,R,A_0) = -1.$$
\end{Theorem} \quad \quad $\square$
\begin{rmrk}
It is easy to see that the ``reversals'' of the above theorems are also true.
\end{rmrk}
Similar problems in homogeneous Thurston geometries
represent huge class of open mathematical problems. For
$\SXR$, $\HXR$, $\NIL$, $\SOL$, $\SLR$ geometries only very few results are known
(e.g. \cite{CsSz16}, \cite{Cs23}, \cite{MSzV},  \cite{MSz}, \cite{MSz12}, \cite{JJS}, \cite{Sz10}, \cite{Sz14-1}, \cite{Sz14-2}, \cite{Sz19}, \cite{Sz24}  ). 
Detailed studies are the objective of
ongoing research.
\medbreak

\end{document}